
\hsize 30pc \vsize 47pc
\input amstex
\input amsppt.sty
\redefine\o{\circ}
\redefine\i{^{-1}}
\define\x{\times}
\define\id{\operatorname{id}}

\def\si{\sigma}
\def\ga{\gamma}

\NoBlackBoxes

\topmatter
\title  An algebraic cell decomposition of the nonnegative part of a flag variety
\endtitle\rightheadtext
{Cell decomposition of $\Cal B_{\ge 0}$}

\author Konstanze Rietsch \endauthor

\address K\. Rietsch, Department of Mathematics, M.I.T., Cambridge, MA 02139
\endaddress

\date August/September 1997\enddate


\abstract
We study the nonnegative part $\Cal B_{\ge 0}$ of the flag variety
$\Cal B$ of a reductive algebraic group, as defined by Lusztig.
Using positivity properties of the canonical basis it is shown that
$\Cal B_{\ge 0}$ has an algebraic cell decomposition indexed by pairs of
elements $w\le w'$ of the Weyl group. This result was conjectured
by Lusztig in \cite{Lu}.
\endabstract

\thanks
The author was supported by an Alfred P. Sloan Dissertation
Fellowship.
\endthanks

\endtopmatter
\document

\head Introduction \endhead

The theory of total positivity for reductive algebraic groups $G$
was introduced by Lusztig in the fundamental paper \cite{Lu}. While
the definitions of the totally positive parts of $G$ and other
related varieties can be stated in elementary ways, many of their
properties are only proved using the deep positivity properties of
canonical bases. The cell decomposition of the totally nonnegative
part $\Cal B_{\ge 0}={(\Cal B_{>0})^{}}^-$ of the flag variety
proved here has been conjectured in \cite{Lu} and is another
example of this phenomenon. The idea behind this cell decomposition
is the following. By Proposition 8.12 in \cite{Lu}, $\Cal B_{>0}$
is a connected component of the intersection of two opposed big
cells. Here we generalize this result by identifying a `totally
positive' connected component in any intersection of two opposed
Bruhat cells and showing that it is topologically a cell. These
connected components are then the cells in the proposed cell
decomposition of $\Cal B_{\ge 0}$. In particular, $\Cal B_{>0}$
becomes the (unique) open cell, and the $0$-dimensional cells are
simply $\{\dot w\cdot B^+\}$ for Weyl group elements $w$. For
general intersections of opposed Bruhat cells it is not obvious
that a positive part exists, that is, that their intersection with
$\Cal B_{\ge 0}$ is nonempty. Showing this is in some sense the
heart of the proof, and is where we require the use of canonical
bases (see Lemma 6).

The parameterization of this decomposition comes from indexing the
Bruhat cells and opposite Bruhat cells by Weyl group elements in an
appropriate way. Then pairs of opposed Bruhat cells with nonempty
intersection (of dimension $\ell(w')-\ell(w)$) are labelled
precisely by pairs $(w,w')$ with $w\le w'$ (see \cite{Lu2}).

\head Preliminaries \endhead
Let $G$ be a reductive linear algebraic group split over $\Bbb R$
with fixed pinning $(T,B^+,B^-,x_i,y_i;i\in I)$ (see \cite{Lu}). In
the following all varieties will be identified with their $\Bbb
R$-valued points. Let $\Cal B$ be the variety of all Borel
subgroups of $G$, and $W=N_G(T)$ denote the Weyl group with longest
element $w_0\in W$. Let $\{s_i|\ i\in I\}$ be the set of simple
reflections in $W$ corresponding to the pinning. Write $\ell(w)$
for the length of $w\in W$, and $\dot w$ for a representative of
$w$ in $N_G(T)$. We start by recalling some results of Lusztig's
from \cite{Lu}.

\proclaim{Total positivity in $U^-$}\rm Let $U^-$ and $U^+$
denote the unipotent radicals of $B^-$ respectively $B^+$. The
totally nonnegative part $U^{-}_{\ge 0}$ of $U^-$ is the semigroup
generated by the set $\{y_i(a)|\ a\in\Bbb R_{\ge 0},\ i\in I\}$.
Analogously, $U^+_{\ge 0}$ is the semigroup generated by
$\{x_i(a)|\ a\in\Bbb R_{\ge 0},\ i\in I\}$.

For every $w\in W$ with reduced expression $w=s_{i_1}\cdots
s_{i_k}$ consider the morphism
$$
\align
 (\Bbb {R^*})^k\quad &\rightarrow \quad U^- \\
 (a_1,\dotsc, a_k)&\mapsto y_{i_1}(a_1)\cdots y_{i_k}(a_k). \\
\endalign
$$
depending on $\bold i:=(i_1,\dotsc, i_k)\in I^k$. The image of the
connected component $\Bbb R_{>0}^k$ under this map clearly lies in
$U^-_{\ge 0}$. We denote it by $U^-(w)$. By \cite {Lu, 2.7} it is
independent of the choice of reduced expression $\bold i$. In fact,
$U^-(w)$ can also be characterised as
$$
U^-(w)=B^+ \dot w B^+ \cap U^-_{\ge 0}.
$$
In particular, the $U^-(w)$ form a cell decomposition of $U^-_{\ge
0}$.
$$
U^-_{\ge 0}=\bigsqcup_{w\in W} U^-(w).
$$
For $w=w_0$, the component $U^-(w_0)$ is itself a semigroup. It is
called the `totally positive submonoid' $U^-_{>0}$ of $U^-$. The
closure of $U^-_{>0}$ turns out to coincide with $U^{-}_{\ge 0}$,
(see \cite {Lu, 4.2}).

In all of the above, $U^-$ can be replaced by $U^+$, with the
$y_i$'s replaced by the $x_i$'s. In this way $U^+_{>0}$, $U^+_{\ge
0}$ and the $U^+(w)$'s are defined.
\endproclaim

\proclaim{Total positivity in $\Cal B$}\rm
We denote the conjugation action of $g\in G$ on $\Cal B$ by
$$
\align
 g:&\Cal B\longrightarrow \Cal B\\
 &B\mapsto g\cdot B:=g B g\i.\\
\endalign
$$
The totally positive part $\Cal B_{>0}$ of the flag variety is then
defined as the orbit $U^-_{>0}\cdot B^+$ of $B^+\in \Cal B$ under
conjugation by elements of $U^-_{>0}$. And this definition turns
out to be symmetric in $U^+$ and $U^-$ (see \cite {Lu, 8.7}),
$$
\Cal B_{>0}=U^-_{>0}\cdot B^+=U^+_{>0}\cdot B^-.
$$
$\Cal B_{\ge 0}$ is by definition the closure of $\Cal B_{>0}$ in
$\Cal B$.

The proposed cell decomposition of $\Cal B_{\ge 0}$ comes from the
two opposed Bruhat decompositions of $\Cal B$. We first set up some
notation. Write
$$
B\overset w\to\longrightarrow B'\quad\ \text{for}\ \ B,B'\in\Cal B\
\ \text{and}\ \ w\in W
$$
if $(B,B')$ is conjugate under $G$ (acting diagonally) to the pair
$(B^+,\dot w\cdot B^+)$. Let $w,w'\in W$. We use the following
notation for Bruhat cells with respect to $B^-$ and $B^+$,
$$
\alignat{2}
 &\Cal C^-_{w'}&:=\{B|B^-\overset{w'}\to\longrightarrow B\}&\\
 &\Cal C^+_{w}&\ :=\{B|B^+\overset{w_0 w}\to\longrightarrow B\}.&\\
\endalignat
$$
Note that $\Cal C^+_{w}$ has codimension $\ell(w)$ in $\Cal B$. The
intersection $\Cal C^-_{w'}\cap\Cal C^+_{w}$ is nonempty precisely
if $w\le w'$, in which case it is a smooth variety of dimension
$\ell(w')-\ell(w)$. This was shown by Lusztig in \cite {Lu2}, and
also by Dale Peterson, \cite{P}. We write
$$
\Cal R_{w,w'}=\Cal C^-_{w'}\cap\Cal C^+_w.
$$

Our aim is to prove that $\Cal R_{w,w';>0}:=\Cal R_{w,w'}\cap\Cal
B_{\ge 0}$ is a cell of dimension $\ell(w')-\ell(w)$. More
precisely we propose that it is homeomorphic to $\Bbb
R_{>0}^{\ell(w')-\ell(w)}$ by a homeomorphism that extends to a
(real) algebraic morphism $(\Bbb R^*)^{\ell(w')-\ell(w)}\rightarrow
\Cal R_{w,w'}$. This was conjectured by Lusztig
(see \cite{Lu, 8.15} and also \cite{Lu2}).
\endproclaim

\proclaim {Key Example}\rm
If $w=1$ or if $w'=w_0$ in $\Cal R_{w,w'}$ then the conjecture
follows directly from Lusztig's cell decomposition of $U^+_{\ge
0}$, respectively $U^-_{\ge 0}$. In that case
$$
\align
 &\Cal R_{1,w';>0}=U^+(w_0 w' w_0)\cdot B^-\\
 &\Cal R_{w,w_0;>0}=U^-(w_0w)\cdot B^+.
\endalign
$$
The desired algebraic maps $(\Bbb R^*)^{\ell(w')}\rightarrow\Cal
R_{1,w'}$ and $(\Bbb R^*)^{\ell(w_0)-\ell (w)}\rightarrow\Cal
R_{w,w_0}$ are given by
$$
\align
 (a_1,\dotsc,a_k)&\mapsto x_{i_1}(a_1)\dotsc x_{i_k}(a_k)\cdot B^-\\
 (b_1,\dotsc,b_m)&\mapsto y_{j_1}(b_1)\dotsc y_{j_m}(b_m)\cdot B^+\\
\endalign
$$
where $s_{i_1}\dotsc s_{i_k}=w_0 w' w_0$ and $s_{j_1}\dotsc
s_{j_m}=w_0 w$ are reduced expressions.
\endproclaim

\head
Proof of the Cell Decomposition
\endhead

The proof follows the recursive procedure for studying $\Cal
R_{w,w'}$ of Kazhdan and Lusztig, \cite{K-L}, examining what
happens to the positive part at every step. For $B\in\Cal B$, we
will sometimes write $B\ge 0$ instead of $B\in\Cal B_{\ge 0}$.

Let $w, v\in W$ with $\ell(wv)=\ell(w)+\ell(v)$. Denote by
$$
\phi_{w, v}:\ \ \Cal C^-_{wv}
       \longrightarrow \Cal C^-_{w} \qquad\text{and}\qquad
\phi^{w,v}:\ \ \Cal C^+_{w}\longrightarrow \Cal C^+_{wv}
$$
the algebraic maps defined by
$$
\align
 &B^-\overset{w}\to\longrightarrow \phi_{w,v}(B)
 \overset{v}\to\longrightarrow B\\
 &B^+\overset{w_0 w v}\to\longrightarrow \phi^{w,v}(B)\overset{v\i}
 \to\longrightarrow B.
\endalign
$$
We start by computing $\phi_{w,v}$ in the key example.

\proclaim{1. Lemma}
Let $w'\in W$, $\ B\in\Cal R_{1,w',>0}$, and $v\in W$ such that
$\ell(w'v)=\ell(w')-\ell(v)$. Then $\phi_{w'v,v}(B)\in\Cal
R_{1,w'v;>0}$.
\endproclaim

\demo{Proof}
By the key example, $B=x\cdot B^-$ for some $x\in U^+(w_0 w' w_0)$.
We write $x$ as $x_{i_1}(a_1)\dotsc x_{i_m}(a_m)$ where
$a_1,\dotsc, a_m\in\Bbb R_{>0}$ and $s_{i_1}\dotsc s_{i_m}=w_0 w'
w_0$ is a reduced expression with $s_{i_{k+1}}\cdots s_{i_m}=w_0 v
w_0$. Then $\phi_{w'v,v}(B)= x_{i_1}(a_1)\dotsc
x_{i_{k}}(a_{k})\cdot B^-$ and therefore lies in $U^+(w_0 w'v
w_0)\cdot B^-=\Cal R_{1\le w'v,>0}$. \qed
\enddemo

To try to reduce the study of $\Cal R_{w,w';>0}$ to the key example
$\Cal R_{1,w';>0}$, we conjugate by an element of $U^-(w_0 w\i
w_0)$. Some properties of the resulting map are listed in Lemma 2.

\proclaim{2. Lemma}
Let $w\in W$ and $y\in U^-(w_0 w\i w_0)$.
\roster
\item
If $B\in\Cal B_{\ge 0}$ then $y\cdot B\in\Cal B_{\ge 0}$.
\item
Conjugation by $y$ induces embeddings
$$
\Cal C^+_{w}\rightarrow \Cal C^+_1,\quad
\Cal C^-_{w'}\rightarrow\Cal C^-_{w'},\quad\text{and}\quad
\Cal R_{w,w'}\rightarrow \Cal R_{1,w'}.
$$
\item
For $B\in\Cal C^-_{w'v}$ with $\ell(vw')=\ell(v)+\ell(w')$, we have
$
\phi_{w',v}(y\cdot B)=y\cdot\phi_{w',v}(B).
$
\endroster
\endproclaim

\demo{Proof}
\roster
\item
Conjugation by $y$ preserves $U^-_{\ge 0}\cdot B^+$ therefore also
its closure, which is $\Cal B_{\ge 0}$.
\item
Since $y\in B^+\dot w_0\dot w\i\dot w_0 B^+$ we have
$B^+\overset{w_0 w\i w_0}\to\longrightarrow y\cdot B^+$.
Furthermore, $B\in\Cal C^+_{w}$ implies $y\cdot B^+\overset{w_0
w}\to\longrightarrow y\cdot B$. Therefore
$B^+\overset{w_0}\to\longrightarrow y\cdot B$, since the lengths
add up. The rest follows since $y\in U^-$.
\item
By definition, $B^-\overset
w'\to\longrightarrow\phi_{w',v}(B)\overset v\to\longrightarrow B$.
Conjugation by $y$ gives $B^-\overset w'\to\longrightarrow
y\cdot\phi_{w',v}(B)\overset v\to\longrightarrow y\cdot B$ which
implies (3). \qed
\endroster\
\enddemo

\proclaim {3. Lemma}
Suppose $v,w,w'\in W$, $\ \ell(wv)=\ell(w)+\ell(v)$ and
$\ell(w'v)=\ell(w')+\ell(v)$.

(1) Let $B\in\Cal C^-_{w'v}$. If $B\ge 0$, then $\phi_{w',v}(B)\ge
0$.

(2) Let $B\in\Cal C^+_{v}$. If $B\ge 0$, then $\phi^{w,v}(B)\ge 0$.

(3) If $w\le w'$ then $\phi_{w,v}$ gives rise to an isomorphism
$
\phi:\Cal R_{wv,w'v}\rightarrow \Cal R_{w,w'}
$
that restricts to a bijection
$$
\phi_{>0}:\Cal R_{wv,w'v;>0}\rightarrow \Cal R_{w,w';>0}.
$$

\endproclaim

\demo{Proof}
Note that statement (1) is true for $B\in \Cal R_{1,w'v}\subset\Cal
C^-_{w'v}$, by Lemma 1. Now suppose $B\in\Cal R_{w,w'v;>0}$ where
$w\ne 1$. In that case consider the curve $\Bbb R_{>0}\rightarrow
U^-(w_0 w\i w_0)$ given by
$$
t\mapsto y(t):=y_{j_1}(t) y_{j_2}(t)\dots y_{j_l}(t)
$$
for some reduced expression $s_{j_1}s_{j_2}\cdots s_{j_l}$ of $w_0
w\i w_0$. By Lemma 2 we have $y(t)\cdot B\in\Cal R_{1,w'v;>0}$.
Hence $\phi_{w',v}(y(t)\cdot B)\ge 0$ and by continuity as $t$ goes
to $0$, also $\phi_{w',v}(B)\ge 0$. This implies (1), since $\Cal
C^-_{w'v}=\bigsqcup_{w}\Cal R_{w,w'v}$. (2) follows by symmetry.
(3) is an immediate consequence of (1) and (2), since the inverse
of $\phi$ is given by $\phi^{w,v}|_{\Cal R_{w,w'}}$.
\qed
\enddemo

Let $w\le w'\in W$ and $s$ be a simple reflection such that $w\le
ws$ and $w's\le w'$. By properties of the Bruhat decomposition,
$\phi_{w's,s}$ restricts to
$$
 \pi=\pi_{w,w',s}:\Cal R_{w,w'}\rightarrow \Cal R_{w,w's}\sqcup \Cal R_{w
 s,w's}.
$$
The following lemmas study the behaviour of $\pi_{w,w',s}$ (denoted
$\pi$ if the $w,w'$ and $s$ are clear from context) with respect to
total positivity.

\proclaim{4. Lemma}
For $w,w',s\in W$ as above, the map $\pi=\pi_{w,w',s}$ restricts to
a map
$$
\pi_{>0}:\Cal R_{w,w',>0}\rightarrow \Cal R_{w,w's,>0}.
$$
\endproclaim

\demo{Proof}
Suppose $B$ is an element of $\Cal R_{w,w',>0}$. By Lemma 3.(1) we
know that $B':=\pi(B)$ lies in $\Cal B_{\ge 0}$. It remains to show
that $B'$ lies in $\Cal R_{w,w's}$, as opposed to in $\Cal
R_{ws,w's}$. Choose some $y\in U^-(w_0w\i w_0)$. Then by Lemma 2 we
have $y\cdot B'= y\cdot\phi_{w's,s}(B)=\phi_{w's,s} (y\cdot B)$ and
$y\cdot B\in\Cal R_{1,w';>0}$. Therefore, by Lemma 1, $y\cdot
B'\in\Cal R_{1,w's;>0}$. But this can only be the case if
$B'\in\Cal R_{w,w's}$. If $B'$ were in $\Cal R_{ws,w's}$ then we
would have
$$
 B^+\overset{w_0 w\i w_0}\to\longrightarrow y\cdot B^+
 \overset {w_0 w s}\to \longrightarrow y\cdot B'.
$$
and $y\cdot B'\in\Cal R_{s,w's}$, which is a contradiction.
\qed
\enddemo

\proclaim{5. Total positivity and canonical bases}\rm
We have not yet shown that $\Cal R_{w,w';>0}$ is nonempty. Our
proof of this fact requires the deep positivity properties of
Lusztig's canonical basis. Let $\rho$ be the sum of all fundamental
weights, $V$ an irreducible representation of $G$ with highest
weight $\rho$ and $\eta\in V$ a highest weight vector. Then
Lusztig's canonical basis of the (quantized) universal enveloping
algebra of $U^-$ (with respect to the chosen pinning, or
equivalently the corresponding set of Chevalley generators) gives
rise to a basis $\Bbb B$ of $V$ uniquely determined by the choice
of $\eta$ (see \cite {Lu3}). Using this basis, $\Cal B_{\ge 0}$ can
be characterized as follows.
\endproclaim

\remark{Theorem(Lusztig \cite{Lu, 8.17})} Let $G$ be of simply
laced type and $B\in\Cal B$. Then $B$ lies in $\Cal B_{\ge 0}$ if
and only if the unique line in $V$ stabilized by $B$ is spanned by
a vector $v\in V_{\ge 0}=\sum_{b\in\Bbb B}\Bbb R_{\ge 0}\ b$.
\endremark

\vskip .25 cm
We apply this Theorem in the proof of the next lemma to explicitly
construct elements in $\Cal R_{w,w';>0}$.

\proclaim{6. Lemma}
Let $y\in U^-(w_0 w\i w_0)$ and $B\in\Cal C^+_w$. If $y\cdot
B\in\Cal B_{\ge 0}$, then $B\in\Cal B_{\ge 0}$.
\endproclaim

\demo{Proof}
Note that it suffices to prove this lemma for simply laced groups.
The non-simply laced case then follows by standard arguments (see
\cite {Lu,8.8}). So we assume that $G$ is simply laced.

$\Cal C^+_{w}$ is decomposed as follows.
$$
\Cal C^+_{w}=\bigsqcup_{w'\ge w}{\Cal R_{w,w'}}
$$
We prove the lemma stratum by stratum by induction on
$\ell(w')-\ell(w)$ (for all $w\in W$ simultaneously). In the
starting case of the induction a stratum consists of just a single
element, $\Cal R_{w,w}=\{\dot w_0
\dot w\i\cdot B^+\}$. By \cite{Lu, 8.13}, this element lies in $\Cal
B_{\ge 0}$, thus the statement of the lemma holds.

Consider $\Cal R_{w,w'}$ with $w<w'$. First we show that we can
assume the existence of a simple reflection $s$ such that $w\le ws$
and $w' s\le w'$. Let $v\in W$ be maximal such that
$\ell(wv)=\ell(w)+\ell(v)$ and $\ell(w'v)=\ell(w')+\ell(v)$. We
reduce the case $B\in\Cal R_{w,w'}$ to $B'\in\Cal R_{wv,w'v'}$. Let
$B\in\Cal R_{w,w'}$ and $y\in U^-(w_0 w\i w_0)$ such that $y\cdot
B\ge 0$. Then $B':=\phi_{w',v}(B)\in\Cal R_{wv,w'v}$ and $y\cdot
B'=\phi_{w',v}(y\cdot B)\ge 0$, by Lemmas 2 and 3. Choose $y'\in
U^-(w_0 v\i w_0)$, so that $y'y\in U^-(w_0 v\i w\i w_0)$. Then also
$y'y\cdot B'\ge 0$. Thus if the lemma holds for $y'y\in U^-(w_0 v\i
w\i w_0)$ and $B'\in\Cal R_{wv,w'v}$, then by Lemma 3.(3) it also
holds for $y\in U^-(w_0 w\i w_0)$ and $B\in\Cal R_{w,w'}$.
Therefore we can replace $w,w'$ by $wv,w'v$. By construction, any
simple reflection that increases the length of this new $w$ must
decrease the length of the new $w'$. And such a simple reflection
exists, since $w<w'\le w_0$.

So we are reduced to considering strata $\Cal R_{w,w'}\subset\Cal
C^+_{w}$ for $w,w'$ with a simple reflection $s$ such that $w\le
ws$ and $w's\le w'$. This is the case in which we can apply the
induction hypothesis. Consider the commutative diagram
$$
\alignat{2}
 &\qquad\qquad\Cal R_{w,w'}
                    &\overset y\to\longrightarrow &\quad\Cal R_{1,w}\\
 &\qquad\qquad\pi\downarrow &                  &\quad\downarrow \pi_1\\
 &\Cal R_{w,w' s}\sqcup\Cal R_{ws,w' s}&\overset y\to\longrightarrow &
                           \quad\Cal R_{1,w' s}\sqcup\Cal R_{s,w's}.\\
\endalignat
$$
where $\pi$ and $\pi_1$ are restrictions of $\phi_{w's,s}$ as
before, and the horizontal maps refer to conjugation by $y$.
Suppose $B\in\Cal R_{w,w'}$ and $y\in U^-(w_0 w\i w_0)$ such that
$y\cdot B\ge 0$. Then $y\cdot
\pi(B)=\pi_1(y\cdot B)$ and lies in $\Cal R_{1,w's;>0}$, by Lemma 1.
Thus, by the induction hypothesis, $\pi(B)\in\Cal R_{w,w's;>0}$
(and since $y$ takes $\Cal R_{ws,w's}$ to $\Cal R_{s,w's}$).

Consider the fibers $F=\pi\i(\pi(B))$ and
$F_1=\pi_1\i(y\cdot\pi(B))$ and their nonnegative parts
$F_{>0}:=F\cap \Cal B_{\ge 0}$ and $F_{1;>0}:=F_1\cap \Cal B_{\ge
0}$. Then $B$ is an element of $F$ such that $y\cdot B\in
F_{1;>0}$. It remains to show that $B\ge 0$. This follows from the
following claim.

\remark{Claim}
Let $B'\in\Cal R_{w,w's;>0}$, and let $F$ and $F_1$ be the fibers
$\pi\i(B')$, respectively $\pi_1\i(y\cdot B')$. The isomorphism
$y:F\rightarrow F_1$ given by conjugation with $y$ restricts to a
bijection $F_{>0}\rightarrow F_{1;>0}$ of the positive parts.
\endremark

\vskip .25cm
The proof of this claim uses the characterization of $\Cal B_{\ge
0}$ in terms of canonical bases to construct $F_{>0}$. Let $g\in
B^+\dot w_0\dot w T$ such that $g\cdot B^+=B'$. We apply $g$ to the
highest weight $\eta$ of $\rho$-representation. Since $B'\ge 0$ we
can assume that, by Lusztig's Theorem,
$$
g.\eta\in\sum_{b\in\Bbb B}\Bbb R_{\ge 0}\ b= V_{\ge 0}.
$$
Now in explicit terms $F=\{g y_i(a)\cdot B^+|\ a\in\Bbb
R\setminus\{0\}\ \}$ where $s=s_i$. To determine which elements of
$F$ lie in $F_{>0}$ we need to study the $g y_i(a).\eta$ for
$a\in\Bbb R^*$. Let $f_i$ be the Chevalley generator of the Lie
algebra of $U^-$ such that $y_i(a)=\exp(a f_i)$. Then in the
$\rho$-representation $y_i(a).\eta=\eta + a f_i.\eta$ and
$f_i.\eta\in\Bbb B$ is the unique canonical basis element in the
$(s_i.\rho)$-weight space. We have therefore $g. (f_i.\eta)=g \dot
s_i.\eta$ for a suitable choice of $\dot s_i$. It is easily
verified that $g\dot s_i\cdot B^+=\phi^{w,s}(g\cdot B^+)$ and thus
lies in $\Cal B_{\ge 0}$. So $g. (f_i.\eta)\in \si V_{\ge 0}$ for
some sign $\si\in\{\pm 1\}$. Therefore
$$
g\ y_i(a).\eta=g.\eta+ a\ g.(f_i.\eta) \in V_{\ge 0}
\qquad\text{whenever}\ \ \si a\in\Bbb R_{>0},
$$
and the connected component
$$
F^0:=\{g\ y_i(a)\cdot B^+|\ \si a\in\Bbb R_{>0}\}\subseteq F_{>0}
$$
of $F$ lies in $F_{>0}$. Now $F_{1;>0}$ is indeed a connected
component of $F_1$, as can easily be computed using the key
example. So its image under $y\i$ must in fact coincide with $F^0$,
and we get inclusions
$$
y\cdot F_{>0}\subseteq F_{1;>0}=y\cdot F^0\subseteq y\cdot F_{> 0}.
$$
Thus $F_{>0}=F^0$, and conjugation by $y$ induces a bijection
$F_{>0}\rightarrow F_{1;>0}$.
\qed
\enddemo

\proclaim{7. Lemma} Let $w<w'\in W$ and $s$ a simple
reflection such that $w\le ws$ and $w' s\le w'$. There exists a
real algebraic map $\psi:\Cal R_{w,{w's}}\x\Bbb R^*\rightarrow\Cal
R_{w,w'}$ that restricts to a homeomorphism
$$
\psi_{>0}:\Cal R_{w,{w's};>0}\x\Bbb
R_{>0}\rightarrow\Cal R_{w,w';>0}
$$
and such that the following diagram commutes
$$
\alignat{2}
 &\Cal R_{w,{w's}}\x\Bbb R^*&\overset{\psi}\to\longrightarrow & \Cal
 R_{w,w'}   \\
 & \qquad pr_1 \downarrow & &\downarrow \pi \\
 &\qquad \Cal R_{w,w's}&\overset \iota\to\longrightarrow &\Cal R_{w,w's}\sqcup\Cal
 R_{ws,w's}.
\tag *\endalignat
$$
Here $pr_1$ is projection onto the first factor,
$\pi=\pi_{w,w',s}$, and $\iota$ is the obvious inclusion.
\endproclaim
\demo{Proof}
First note that this statement is true for $w=1$ by the key
example. Here $\psi=\psi_1:\Cal R_{1,w's}\x\Bbb R^*\rightarrow\Cal
R_{1,w'}$ is given by $(x\cdot B^-,a)\mapsto xx_i(a)\cdot B^-$,
where $x\in U^+\cap B^-\dot w_0\dot w'\dot s\dot w_0 B^-$ and
$s_i=w_0 s w_0$. It clearly has all the properties required in the
lemma.

Let $y\in U^-(w_0 w\i w_0)$. We consider the following commutative
diagram, where $\pi_1=\pi_{1,w,w'}$ and $y$ stands for conjugation
by $y$.
$$
\alignat{2}
 &\Cal R_{w,w's}\x\Bbb R^*\ \ &\qquad\quad &\Cal R_{w,w'} \\
 &y\x\id\downarrow          &\overset {\tilde\psi}\to
 \searrow\quad\  &\downarrow y     \\
 &\Cal R_{1,w's}\x\Bbb R^* &\overset{\psi_1}\to
                          \longrightarrow \quad&\Cal R_{1,w'}\\
 &\quad pr_1\downarrow        &            &\downarrow\pi_1   \\
 &\qquad\Cal R_{1,w' s}&\overset\iota\to\longrightarrow\quad&\Cal R_{1,w's}\sqcup
  \Cal R_{s,w's}
\tag **\endalignat
$$
By the lower half of this diagram we see that the image of
$\tilde\psi:=\psi_1\o(y\x\id)$ lies in $\pi_1\i(y\cdot\Cal R_{w,
w's})=y\cdot(\pi\i(\Cal R_{w,w's}))\subseteq y\cdot(\Cal
R_{w,w'})$. Therefore
$
\psi(B,a):= y\i\cdot\psi_1(y\cdot B,a)
$
defines a real algebraic map $\psi:\Cal R_{w,w's}\x\Bbb
R^*\rightarrow\Cal R_{w,w'}$ with the property
$\tilde\psi=y\o\psi$.

Next we study the restriction of $\psi$ to $\Cal R_{w,
w's;>0}\x\Bbb R_{>0}$. Consider an element $(B,a)\in\Cal R_{w,
w's;>0}\x\Bbb R_{>0}$. Then $\psi(B,a)\in\Cal R_{w,w'}$ and, by the
properties of $\psi_1$, $y\cdot\psi(B,a)=\psi_1(y\cdot B,a)$ lies
in $\Cal R_{1,w';>0}$. Therefore $\psi(B,a)\ge 0$, by Lemma 6. Thus
the restriction of $\psi$ gives rise to a continuous map
$\psi_{>0}:\Cal R_{w,w's;>0}\x\Bbb R_{>0}\rightarrow \Cal
R_{w,w';>0}$. Its inverse should be given by $B'\mapsto
(y\i\x\id)(\psi_1\i(y\cdot B'))$. The first component of this map
is just $y\i\o\pi_1\o y=\pi$ and the second component equals to the
second component of $\psi_1\i\o y$. Thus by Lemmas 4 and 2.(1),
$\psi_{>0}\i(B'):=(y\i\x\id)(\psi_1\i(y\cdot B'))$ lies in $\Cal
R_{w,w's;>0}\x\Bbb R_{>0}$ (for $B'\in\Cal R_{w,w';>0}$), and
$\psi_{>0}$ is a homeomorphism.

It remains to note that the diagram (*) commutes. This follows
since (*) can be obtained from the lower half of the commutative
diagram (**) by conjugating all the maps by $y\i$ and restricting.
\qed
\enddemo

\proclaim{8. Proposition}
Let $w<w'\in W$ and $m:=\ell(w')-\ell(w)$. Then there exists a real
algebraic map $\ga:(\Bbb R^*)^m\rightarrow\Cal R_{w,w'}$ such that
its restriction defines a homoemorphism $\ga_{>0}:\Bbb
R_{>0}^m\rightarrow
\Cal R_{w,w';>0}$.
\endproclaim

\demo{Proof}
Let $v\in W$ be of maximal length such that
$\ell(wv)=\ell(w)+\ell(v)$ and $\ell(w'v)=\ell(w')+\ell(v)$. By
Lemma 2 it suffices to prove the Proposition for $\Cal R_{wv,w'v}$.
Thus it suffices to consider the case where we have a simple
reflection $s$ such that $w\le ws$ and $w's\le w'$.

The Proposition is proved by induction on $\ell(w')-\ell(w)$. For
the start of induction we can assume that $w'=ws$ for some simple
reflection $s$. In this case Lemma 7 applies and (since $\Cal
R_{w,w's}=\Cal R_{w,w}$ is a single point) gives the desired map
$\Bbb R^*\rightarrow\Cal R_{w,w s}$.

Now for general $w<w'$ and $s$ with $w\le ws$ and $w' s\le w'$ we
use the algebraic map $\ga': (\Bbb R^*)^{m-1}\rightarrow
\Cal R_{w,w's}$ given by the induction hypothesis.
The map $\ga: (\Bbb R^*)^m\rightarrow\Cal R_{w,w'}$ is then defined
as the composition
$$
\ga:(\Bbb R^*)^m\cong(\Bbb R^*)^{m-1}\x\Bbb R^* \overset \ga'\x\id \to\longrightarrow
     \Cal R_{w,w's}\x\Bbb R^*\overset\psi\to\longrightarrow
     \Cal R_{w,w'},
$$
with $\psi$ as in Lemma 7. It is clearly algebraic and restricts to
a homeomorphism $\ga_{>0}=\psi_{>0}\o(\ga'_{>0}\x\id):\Bbb
R_{>0}^m\cong\Bbb R^{m-1}_{>0}\x\Bbb R_{>0}\rightarrow\Cal
R_{w,w';>0}$.
\qed
\enddemo

\remark{Remark}
It remains to note that all $W$-conjugates of $B^+$ lie in $\Cal
B_{\ge 0}$, by \cite{Lu, 8.13}. Therefore $\Cal R_{w,w;>0}=\Cal
R_{w,w}=\{\dot w_0 \dot w\cdot B^+\}$. With this, the proof of the
cell decomposition is complete.
\endremark

\remark{Acknowledgements}
The author would like to thank George Lusztig for suggesting the
problem.
\endremark

\enddocument